\documentclass[12pt]{article}
\usepackage{amsmath,amsfonts,amssymb,amscd,url}

\newcommand{\HH}{{\mathbb H}}

\newcommand{\kk}{{\bf k}}

\newtheorem{theorem}{Theorem}
\newtheorem{lemma}{Lemma}

\title{Two-sided fundamental theorem of affine geometry}
\author{A.\,G.\,Gorinov}
\date{}
\usepackage{fullpage} 
\begin{document}
\maketitle

\begin{abstract}
The fundamental theorem of affine geometry says that a self-bijection $f$ of a finite-dimensional affine space over a possibly skew field takes left affine subspaces to left affine subspaces of the same dimension, then $f$ of the expected type, namely $f$ is a composition of an affine map and an automorphism of the field. We prove a two-sided analogue of this: namely, we consider self-bijections as above which take affine subspaces  affine subspaces  but which are allowed to  take left subspaces to right ones and vice versa. We show that these maps again are of the expected type.
\end{abstract}

\section{Introduction}

Let $\kk$ be an associative division algebra. The ``fundamental theorem of affine geometry'' is the statement that if a bijection $f:\kk^n\to\kk^n$ takes every left affine subspace to a left affine subspace of the same dimension, then $f$ is an affine map, composed with a map induced by an automorphism of $\kk$, provided $\infty>n\geq 2$, see e.g.\,\cite[Chapter 2]{artin}  Here we prove a two-sided version of this. We will say that an affine subspace $\subset\kk^n$ is {\it purely left} (respectively {\it purely right}) iff it is left but not right (respectively right but not left). Affine subspaces which are both right and left will be called {\it two-sided}.

\begin{theorem}\label{main}
Let $f:\kk^n\to\kk^n$ be an injective map which takes every left or right affine subspace to a left or right affine subspace. Then the following holds:

(1) The map $f$ is bijective and preserves dimensions.

Now let $P$ be a two-sided plane.

(2) Suppose $n\geq 3$. If there exists a purely left line $L\subset P$ such that $f(L)$ is left, then the image of every left affine space is a left affine space of the same dimension. The map $f$ can be written as $f(x)=\sigma( g(xa)+b)$ where $a\in \kk,b\in\kk^n$, $\sigma$ is an automorphism of $\kk$ (applied component-wise) and $g:\kk^n\to \kk^n$ is a map of left $\kk$-vector spaces defined over the centre $Z(\kk)$ of $\kk$.  

(3) Suppose $n\geq 2$. If for every purely left line $L\subset P$ the image $f(L)$ is right, then $\kk$ admits an anti-automorphism, which we will denote by $\varepsilon$. The composition of $f$ with $\varepsilon$ applied component-wise is then as in part (2) above.

\end{theorem}

So in particular, injective maps $\kk^n\to \kk^n$ which take affine subspaces to affine subspaces so that at least one left affine subspace is taken to a right one, or vice versa, are possible iff $\kk$ has an anti-automorphism.

\smallskip

In \cite[Problem 2002-10]{arn} V.\,I\,Arnold asked whether the statement is true if $f$ is assumed to be a homeomorphism and $\kk=\HH$. The author was able to provide an affirmative answer \cite{plin}. V.\,I\,Arnold then asked the author whether the statement remains true for arbitrary division algebras. In this note we show that this is indeed the case when $n\geq 3$.

\section{Linear algebra over a division ring}

In order to prove the theorem we need to make a few observations first. The author was unable to find a reference for these, so proofs are provided for completeness.

Let $e_1,\ldots,e_n$ be the standart basis of $\kk^n$.

\begin{lemma}\label{dimleft}
Suppose $V\subset \kk^n$ is a left vector subspace. Let $I$ be the maximum subset of $\{1,\ldots, n\}$ such that $V$ intersected with the (two-sided) subspace spanned by $e_i,i\in I$ is 0. The left dimension of $V$ is $n-\#I$.

The same is true for right vector subspaces.
\end{lemma}

{\bf Proof.} We will consider the case when $V$ is a left vector subspace. Let $k$ be the left dimension of $V$. Choose a basis $x_1,\ldots x_k$ of $V$ and form a matrix $M$ whose rows are $x_1,\ldots x_k$. By using row transformations we can transform $M$ to get a matrix $M'$ such that the $i$-th row of $M'$ has 1 at $j_i$-th place and zeroes before that, where $1\leq j_1<j_2<\ldots <j_k\leq n$ is an increasing sequence of integers. The rows of $M'$ will still be a basis of $V$.

Now set $I=\{1,\ldots, n\}\setminus\{j_1,\ldots, j_k\}$ and set $E_I$ to be the span of $e_i,i\in I$. Clearly, $E_1+V=\kk^n$, so using the exact sequence of left $\kk$-modules
$$0\to E_1\cap V\to E_1\oplus V\to E_1+V\to 0$$ we see that $E_1\cap V=0$ and $k=\dim V=n-\#I$. Using a similar sequence we see that if $E_1'$ is the span of $E_1$ and some $e_j,j\not\in I$, then $E_1'\cap V\neq 0$. $\square$

\begin{lemma}\label{dimleftright}
Suppose $V\subset\kk^n$ is a left vector subspace which contains a right vector subspace $W$. Then the (left) dimension of $V$ is $\geq$ the (right) dimension of $W$. Moreover, the inequality is strict if $W\subsetneq V$.

The same is true with ``left'' and ``right'' interchanged. As a corollary, for two-sided affine or vector subspaces the left and right dimensions coincide.
\end{lemma}

{\bf Proof.} Let $\dim V$ and $\dim W$ be the (left) dimension of $V$ and the (right) dimension of $W$. The inequality $\dim V\geq \dim W$ follows straight from Lemma \ref{dimleft}. To prove that $\dim V\geq \dim W$ once $W\subsetneq V$ let $E_1$ and $I$ be as in the proof of that lemma. Set $E_2$ to be the span of $e_j,j\not\in I$.

The projection $p:\kk^n=E_1\oplus E_2\to E_2$ is a map of $\kk$-bimodules, so when we restrict it to $W$ we get a map of right $\kk$-modules. Moreover, $p$ is injective when restricted to $V$. So if $W\subsetneq V$, then $p(W)$ is a proper right subspace of $E_2\cong \kk^{k}$, which cannot have dimension $\geq k$. $\square$

So from now on we do not distinguish between the left and right dimensions of an affine subspace of $\kk^n$: when both make sense, they coincide.

\medskip

{\bf Proof of part (1) of Theorem \ref{main}.} Suppose $A\subset \kk^n$ is a left affine space. There is a flag
$$A_0\subsetneq A_1\subsetneq A_2\subsetneq \cdots \subsetneq A_n=\kk^n$$ of left affine spaces which contains $A$ and such that $\dim A_i=i$. Let $V_i$ be the vector space associated to $A_i$. We still have strict inclusions
$$0=V_0\subsetneq V_1\subsetneq V_2\subsetneq \cdots \subsetneq V_n=\kk^n.$$

When we apply $f$ to this flag, we get a sequence of vector spaces. These may be left or right, but their dimensions jump at each step by Lemma \ref{dimleftright}. So $\dim f(V_i)=\dim V_i$ for all $i$, and so $\dim f(A)=\dim A$.

This shows that although $f$ might turn left subspaces into right ones and vice versa, it preserves dimensions. In particular, $f(\kk^n)=\kk^n$, which shows that $f$ is surjective. $\square$

\medskip

\begin{lemma}\label{preserve2sided}
An affine subspace $A\subset\kk^n$ is two-sided iff for each affine line $L\subset \kk^n$ (left or right) the intersection $A\cap L$ is $\varnothing$, one point or $L$.
	
As a corollary, if $f$ is as in Theorem \ref{main}, then $f$ takes every two-sided  affine subspace to a two sided affine subspace, and hence takes a purely left or right affine subspace to a purely left or right affine subspace.
\end{lemma}

{\bf Proof.} It suffices to consider the case when $A$ is a vector subspace of $\kk^n$, so let us assume that. Saying that $A$ is left but not right is equivalent to saying that there is a non-zero $x\in A$ such that the right affine line $L$ spanned by $x$ is not contained in $A$. If this is the case, then $L\cap A\neq \varnothing$ a point or $L$. $\square$

\begin{lemma}\label{class2sided}
Suppose $V\subset \kk^n$ is a two-sided vector subspace of dimension $k$. Then there is an isomorphism of $\kk$-bimodules $\kk^n\to \kk^n$ which takes $V$ to $\kk^k$.
\end{lemma}

{\bf Proof.} Let us consider $V$ as a left vector subspace and $E_1, E_2, p:\kk^n\to E_2$ be as in the proofs of Lemmas \ref{dimleft} and \ref{dimleftright}. We have $\kk^n=E_1\oplus V=E_1\oplus E_2$ as $\kk$-bimodules. Moreover, $p$ restricted to $V$ gives an injective map of $\kk$-bimodules $V\to E_2$ of the same dimension, hence an isomorphism. $\square$

\begin{lemma}\label{notmorethanone}
A purely left affine 2-plane $P$ cannot contain more than one right line through any of its points.
\end{lemma}

{\bf Proof.} If $P$ contains two distinct right lines through a point, then $P$ contains a right 2-plane, which by Lemma \ref{dimleftright} implies that $P$ is two-sided. $\square$

\section{Proof of part (2) of Theorem \ref{main}}

Let $f:\kk^n\to\kk^n$ be the map from the theorem.

\begin{lemma}
Let $P$ be a two-sided affine 2-plane and let $L\subset P$ be a purely left affine line such that $f(L)$ is purely left. Let $P'\neq P$ be a two-sided affine 2-plane such that $P'\cap P$ is a line which contains some $X\in L$. Then $f$ takes every purely left affine line $L'\subset P'$ through $X$ to a purely left affine line and ditto for every purely left affine line $L''\subset P$ through $X$.
\end{lemma}

{\bf Proof.} The intersection $P\cap P'$ is two-sided, so $P\cap P'\neq L$. Let $L'\subset P'$ be a purely left line. As a first step, we want to show that $f(L')$ is purely left.

Let $R$ be the left 2-plane spanned by $L$ and $L'$; $R$ is purely left, as if it were two-sided, so would be $L=P\cap R$ and $L'=P'\cap R$. So by Lemma \ref{preserve2sided} $f(R)$ is purely left or purely right. If it were purely right, then so would be $f(R)\cap f(P)$, as $f(P)$ is two-sided by Lemma \ref{preserve2sided}. But $f(R)\cap f(P)=f(R\cap P)=f(L)$, which is assumed purely left. So we conclude that $f(R)$ is purely left, which implies that $f(L')=f(P'\cap R)=f(P')\cap f(R)$ is left, hence purely left.

By symmetry we conclude that every purely left line $L''\subset P$ through $X$ goes to a purely left line. $\square$

Here is a corollary of this lemma:

\begin{lemma}\label{2sidedplanes}
Let $P$ and $P'$ be affine two-sided 2-planes such that $P\cap P'$ is a line. Then if some purely left affine line $L\subset P$ goes to a purely left affine line under $f$, so does every purely left affine line $\subset P\cup P'$.
\end{lemma}

$\square$

\begin{lemma}\label{purelyleftinplanes}
The map $f$ takes every purely left affine line which is $\subset$ some two-sided 2-plane $P'$ to a purely left affine line.
\end{lemma}

This is the part where we need the assumption $n\geq 3$.

{\bf Proof.} Now let $P$ be as in the theorem. By Lemma \ref{2sidedplanes} it suffices to prove that there is a sequence $$P=P_0,P_1,\ldots, P_m=P'$$ of two-sided affine 2-planes in $\kk^n$ such that $P_i\cap P_{i-1}$ is a line for all $i=1,\ldots, m$. Moreover, it suffices to prove this in the case when $P$ is spanned by $e_1,e_2$ by Lemma \ref{class2sided}.

We proceed by induction on $k$ such that the vector 2-plane $V'$ corresponding to $P'$ is included in the 2-sided subspace $\kk^k$ spanned by $e_1,\ldots, e_k$. Let us first consider the case $k=2$. The plane $P'$ is then
$$\{(x,y,a_3,\ldots, a_n)\mid x,y\in\kk\}$$
for some fixed $a_3,\ldots, a_n\in \kk$. The intersection of this subspace with \begin{equation}\label{auxplane}
\{(x',a_2,y',a_4,\ldots, a_n)\mid x',y'\in\kk\}
\end{equation} is a 2-sided line. So both $P'$ and $P''$ given by the same equations as $P'$ but with $a_3$ replaced by 0 intersect (\ref{auxplane}) in a line. Repeating this for all other $a_i$ we construct a sequence of two-sided 2-planes with the required properties.

Suppose now $V'\subset \kk^{k+1}$ but $\not\subset \kk^k$. Then $V'\cap \kk^k$ is a (two-sided) vector 1-subspace $l$. We have $P'=V'+x+ae_{k+1}$ where $a\in\kk$ and $x\in \kk^k$. But $V'+\kk^k=\kk^{k+1}\ni e_{k+1}$ , so $P'=V'+x'$ with $x'\in \kk^k$. This means that $P'\cap \kk^k=l+x'$. Now we can take $P''\subset \kk^k$ to be any two-sides 2-plane that contains $l+x'$ and use the induction hypothesis. $\square$

\smallskip

We now prove that all left lines go to left lines under $f$. Let $L$ be a purely left line. We use induction on the least $k$ such that $L\subset $ an affine two-sided $k$-subspace. The case $k=2$ is Lemma \ref{purelyleftinplanes}.

Suppose $L\subset A$, a two-sided affine $k+1$-subspace, but $L\not\subset$ any affine two-sided $k$-subspace. Choose an $X\in A$ and let $A'$ and $L'$ be a two-sided affine $k$-subspace, respectively, a 2-sided line through $X$; these exist, as $A\cong \kk^{k+1}$ by Lemma \ref{class2sided}. We also know that $L\cap L;=L\cap A'=X$. Let $P$ be the left plane through $L$ and $L'$; note that $P$ is purely left, as $L$ can't be $\subset$ a two-sided 2-plane.

Since $\dim P+\dim A'=\dim A+1$, the intersection $P\cap A'$ is a left line $L''$. This line cannot be two-sided, as $P$, being purely left, cannot contain more than one right line by Lemma \ref{notmorethanone}, and it already contains $L'$.

The image $f(P)$ contains purely left lines $f(L')$ (as $L$ is two-sided) and $f(L'')$ (here we use the induction hypothesis, $L''$ being a subspace of $A$). So $f(P)$ is purely left, again using Lemma \ref{notmorethanone}, and so $f(L)$ is left: $f(P)$ contains $f(L')$ and so it cannot contain more right lines.

\smallskip

So we see that $f$ takes left lines to left lines under the assumption of part (2) of theorem \ref{main}. This then easily implies that all left affine subspaces of dimension $k$ go to left affine subspaces of dimension $k$. We can now use the usual fundamental theorem of affine geometry to conclude that $f$ can be written as $$x\mapsto \sigma(g(x))+v$$ where $\sigma$ is an automorphism of $\kk$ (applied component-wise), $v\in \kk^n$ and $g:\kk^n\to\kk^n$ is a left linear map.

Now we use the fact that $f$ takes right lines to right lines, and hence so does $g$. Let $M$ be the matrix of $g$ in the basis $e_1,\ldots, e_n$. The map $g$ then can be written as $x\mapsto xM$. Here we think of $x$ as a row vector. To prove part (2) of the theorem it remains to show that $M=aM'$, where $a\in \kk$ and all elements of $M'$ are in $Z(\kk)$.

\begin{lemma}\label{2sidedlines}
The element $(x,y)\in\kk^2$ with $x\neq 0$ spans a two-sided vector space iff $x^{-1}y\in Z(\kk)$.
\end{lemma}

$\square$

This lemma shows that every row of $M$ is a constant times a vector in $Z(\kk)^n$. Let these constants be $a_1,\ldots, a_n$. We now show that they are all equal, up to multiplication by an element of $Z(\kk)$.

We have $M=\mathop{\mathrm{diag}}(a_1,\ldots, a_n)N$ where $N$ is an invertible matrix ver $Z(\kk)$. So the map $x\mapsto x \mathop{\mathrm{diag}}(a_1,\ldots, a_n)$ must take right subspaces to right subspaces. Applying this map to $e_i+e_j$, which spans a two-sided line, and using Lemma \ref{2sidedlines} again we see that $a_i^{-1}a_j\in Z(\kk)$ for all $i,j=1,\ldots, n$. This completes the proof of part (2) of the Theorem \ref{main}.

\section{Proof of part (3) of Theorem \ref{main}}

Part (3) of the theorem will follow if we prove this lemma:

\begin{lemma}\label{antiaut}
Let $f:\kk^2\to \kk^2$ be an bijective map which takes affine lines to affine lines. Suppose $f$ of every purely left affine line is purely right one. Then $\kk$ admits an anti-automorphism.
\end{lemma}

{\bf Proof.} As above, we conclude from Lemma \ref{preserve2sided} that $f$ takes two-sided affine lines to two-sided affine lines. So $f$ takes every left affine line to a right one. Also, without restricting generality we may assume that $f(0)=0$.

\begin{lemma}\label{abgroups}
The map $f$ is a homomorphism of abelian groups.
\end{lemma}

{\bf Proof of Lemma \ref{abgroups}.} For $x, y\in \kk^2$ such that $x\neq y$ we will denote the left (respectively right) line through $x$ and $y$ by $L(x,y)$ (respectively $R(x,y)$). We will say that two left lines $\subset\kk^2$ are {\it parallel} iff they do not intersect. Note that if $x,y$ are left linearly independent vectors in $\kk^2$, then the following holds:

\begin{enumerate}
\item $x+y$ is the intersection
point of the left line passing through $x$ and parallel to $y$ and the left line passing through $y$ and parallel to $x$.
\item If $x'\in L(0,x),y'\in L(0,y), x'\neq 0, y'\neq 0, L(x',y)\parallel L(y',x), L(x',y')\parallel L(x,y)$, then $x'=-x,y'=-y$.
\end{enumerate}

There are obvious analogues of these statements for right linearly independent vectors and right lines.

We need to check that then for any $x,y\in\kk^2$
\begin{equation}\label{abgr}
f(x+y)=f(x)+f(y).
\end{equation}

Suppose first that $x,y$ are left linearly independent. Then $f(x),f(y)$ are
right linearly independent: since $f$ is a bijection, the lines $$f(L(0,x))=R(0,f(x)),f(L(0,y))=R((0,f(y))$$ must be different.
Formula (\ref{abgr}) for linearly independent $x,y$ follows from remark 1 above and from the fact that $f$ takes parallel lines to parallel lines.
Similarly, it follows from remark 2 above that if $x,y\in\kk^2$ are linearly independent, then
$f(-x)=-f(x),f(-y)=-f(y)$.

Now suppose that $x,y\in\kk^n$ are non-zero but left linearly dependent, i.e., $y=ax$ for some $a\in\kk,
a\neq 0$. The case $a=-1$ has already been considered, so we can assume $a\neq -1$. Let $z$ be a vector that
is linearly independent with $x$. Since $a\neq -1$, the vectors $x+z,y-z$ are left linearly independent,
so we have $$f(x+y)=f(x+z+y-z)=f(x+z)+f(y-z)=f(x)+f(z)+f(y)+f(-z)=f(x)+f(y).$$ So we have proved \ref{abgr}
for any $x,y\in\kk^2$. $\square$

We now return to the proof of Lemma \ref{antiaut}. Suppose $(x,y)\in\kk^2$ is a non-zero vector. Define a map $\alpha_{x,y}:\kk\to\kk$ by the formula
$f(cx,cy)=f(x,y)\alpha_{x,y}(c)$. All these maps are additive and bijective.
We have $$(f(x,0)+f(0,y))\alpha_{x,y}(c)=f(x,y)\alpha_{x,y}(c)=f(cx,cy)=f(x,0)\alpha_{x,0}(c)+
f(0,y)\alpha_{0,y}(c).$$ Expressing both sides of this equation in terms of $f(1,0)$ and $f(0,1)$ and using the fact that these vectors are right linearly independent (see the proof of the previous lemma) we see that for all $x,y,c\in\kk$ such that $(x,y)\neq (0,0)$ we have $$\alpha_{1,0}(x)\alpha_{x,y}(c)=\alpha_{1,0}(x)\alpha_{x,0}(c)$$ and $$\alpha_{0,1} (y)\alpha_{x,y}(c)=\alpha_{0,1}(y)\alpha_{0,y}(c).$$
This implies that all $\alpha_{x,y}$ coincide. Indeed, if $x\neq 0$ then for any $y\in\kk$
we have $$\alpha_{x,y}=\alpha_{x,0}=\alpha_{x,1}=\alpha_{1,1}.$$ If $y\neq 0$, then
for any $x\in\kk$ we have $$\alpha_{x,y}=\alpha_{0,y}=\alpha_{1,y}=\alpha_{1,1}.$$

Set $\alpha=\alpha_{1,1}$. We now check that $\alpha$ is an anti-automorphism of $\kk$: if $x,y\in\kk$ are both non-zero, then $$f(1,0)\alpha(xy)=f(xy,0)=f(y,0)\alpha(x)=f(1,0)\alpha(y)\alpha(x).$$
$\square$

\begin{flushright}
National Research University\\
Higher School of Economics\\
Moscow, Russia\\
\url{gorinov@mccme.ru, agorinov@hse.ru}
\end{flushright}
\end{document}